\documentclass[12pt]{amsart}
\usepackage{amsfonts}

\usepackage{amscd}
\usepackage{amsmath}

\theoremstyle{plain}
\numberwithin{equation}{section}

\newtheorem{theorem}{Theorem}[section]

\newtheorem{lemma}[theorem]{Lemma}

\theoremstyle{definition}
\newtheorem{definition}[theorem]{Definition}

\newtheorem{remark}[theorem]{Remark}
\newtheorem{example}[theorem]{Example}

\numberwithin{equation}{section}
\newcommand{\oh}{\frac{1}{2}}
\newcommand{\cO}{\mathcal{O}}

\newcommand{\glR}{\mathfrak{gl}(n,\mathbb{R})}
\newcommand{\GL}{\mathrm{GL}}
\newcommand{\GLR}{\mathrm{GL}(n,\mathbb{R})}
\newcommand{\ad}{\mathrm{ad}}
\newcommand{\fa}{\mathfrak{a}}

\newcommand{\fg}{\mathfrak{g}}
\newcommand{\fh}{\mathfrak{h}}
\newcommand{\fk}{\mathfrak{k}}
\newcommand{\fl}{\mathfrak{l}}
\newcommand{\fm}{\mathfrak{m}}
\newcommand{\fn}{\mathfrak{n}}
\newcommand{\fo}{\mathfrak{o}}
\newcommand{\fp}{\mathfrak{p}}
\newcommand{\fq}{\mathfrak{q}}
\newcommand{\fs}{\mathfrak{s}}

\newcommand{\R}{\mathbb{R}}
\newcommand{\C}{\mathbb{C}}

\newcommand{\Z}{\mathbb{Z}}
\newcommand{\Sy}{\mathrm{Sym}}

\theoremstyle{plain}

\begin{document}
\title[Groups and Frames]{Continuous action of Lie groups on
$\mathbb{R}^n$ and Frames}
\author{G. \'{O}lafsson}
\address{Department of Mathematics, Louisiana State University, Baton Rouge, LA
70803, USA}
\email{olafsson@math.lsu.edu}
\thanks{Research supported by NSF grants DMS-0070607 and  DMS-0139783}

\begin{abstract}
Wavelet and frames have become a widely used tool in
mathematics, physics, and
applied science during the last decade. In this article we discuss
the construction of frames for $L^2(\R^n)$ using
the action of closed subgroups $H\subset
\mathrm{GL}(n,\mathbb{R})$ such that $H$ has an open orbit $\cO$
in $\R^n$ under the action $(h,\omega)\mapsto (h^{-1})^T(\omega)$. If
$H$ has the form $ANR$, where $A$ is simply connected and abelian,
$N$ contains a co-compact discrete subgroup and
$R$ is compact containing the stabilizer group of $\omega\in\cO$
then we construct a frame for the space $L^2_{\cO}(\R^n)$ of
$L^2$-functions whose Fourier transform is supported in $\cO$.
We apply this to the case where $H^T=H$ and the
stabilizer is a symmetric subgroup, a case discussed for
the continuous wavelet transform in
\cite{FO02}.
\end{abstract}

\maketitle

\section*{Introduction}

\noindent
The  wavelet transform, and more generally time frequency analysis,
has become a widely used
and studied tool in  mathematics, physics, engineering, and
applied science during the last decade. One of the interesting aspect
is the role played by abstract harmonic analysis and  representation
theory of locally compact groups. In wavelet theory one studies
square integrable representations of semidirect products
$G=\R^n\times_s H$, and in time frequency
analysis representations of the Heisenberg group are used to
understand Gabor frames build from a lattice $\Gamma\subset \mathbb{R}^{2d}$.
In this article we will discuss frames built from the continuous
wavelet transform and discrete subsets $\Gamma$ of $G$.

In the language of representation theory the
continuous wavelet transform on the line is given by
taking
the \textit{matrix coefficients} of the natural representation $\pi$
of the $(ax+b)$-group, i.e., the group of dilations and
translations on the line, on the Hilbert space $L^2(\R)$. Thus
\[\pi (a,b)\psi(x)=|a|^{-1/2}\psi\left(\frac{x-b}{a}\right)=T_bD_a\psi(x)\]
and

\begin{equation}\label{contiwave}
W_\psi (f)(a,b)=
(f\mid \pi (a,b)\psi)=|a|^{-1/2}\int_{\R}f(x)\overline{\psi (\frac{x-b}{a})}\, dx\, .
\end{equation}
Here $T_b:L^2(\R)\to L^2(\R)$ stands for  the unitary isomorphism corresponding
to translation $T_bf(x)=f(x-b)$ and $D_a:L^2(\R)\to L^2(\R)$
is the unitary map corresponding to dilation $D_af(x)=|a|^{-1/2}f(x/a)$, $a\not=0$.

The
discrete wavelet transform is obtained by sampling the wavelet transform,
given by a suitable wavlet $\psi$,
of a function $f$  at points gotten by replacing the full
$(ax+b)$-group by a discrete subset generated by
translation by integers and dilations of the form $a=2^n$:
\begin{eqnarray*}
W_\psi^d(f)(2^{-n},-2^{-n}m) &=&(f\mid \pi ((2^n,m)^{-1})\psi )\\
&=&2^{n/2}\int_{\R}f(x)\overline{\psi (2^n x+m)}\, dx\, .
\end{eqnarray*}
Hence, the corresponding frame is
\begin{equation}\label{frame}
\{\pi ((2^n,m)^{-1})\psi\mid n,m\in \Z\}\, .
\end{equation}
The inverse refers here to the inverse in the $(ax+b)$-group.

This observation, in particular (\ref{contiwave}) is the basis for the generalization of
the continuous and discrete wavelet transform to higher dimensions and more general
settings. For the continuous wavelet transform the relation
to representation theory of the $(ax+b)$-group was already pointed out
by  Grossmann,  Morlet, and Paul
in 1985 \cite{GMP85,GMP86}.
Since then several people have worked on wavelets related to actions
of topological groups acting on $\mathbb{R}^n$. Without trying to be
complete we would like to name the work of  Ali,  Antoine,
and  Gazeau, \cite{AAG91,AAG2000},  Bernier and Taylor \cite
{BT96},  F\"ur and F\"uhr and Mayer \cite{F96,F98,F00,FM01}, and
finally Laugesen, Weaver, Weiss, and Wilson \cite{LWWW2002}. In most
of these cases the group generalizing the $(ax+b)$-group is a semidirect
product $\mathbb{R}^n\times_s H$, where $H$ is a closed subgroup of
$\mathrm{GL}(n,\mathbb{R})$.
For the continuous wavelet transform one often assumes that the group
$H$ has open orbits $\cO_1, \ldots ,\cO_r$ such that
the complement of their union has measure zero. As a further condition for
the existence of wavelet functions, or admissible functions, one
needs that for $\omega$ in an open orbit the stabilizer
$$H^\omega =\{h\in H\mid (h^{-1})^T(\omega)=\omega\}$$
is compact. In \cite{FO02} the case where this condition is not
satisfied were discussed. Instead it was assumed that $H$ is reductive, $H^T=H$,
and that the stabilizer group $H^\omega$ is
a symmetric subgroup, i.e., there exits an involutive automorphism $\tau : H\to H$ such that
$H^\tau_o\subset H^\omega\subset H^\tau$. For any given open orbit $\cO$ we then able to construct a
group $Q=ANR\subset H$ such that the following holds:
\begin{enumerate}
\item There exists points $\omega_0,\ldots ,\omega_k\in \cO$ such that the $Q$-orbits
$\cO_j=Q^T(\omega_j)$, $0\le j\le k$, are open and $\cO\setminus (\cO_1\cup
\ldots \cO_k)$ has measure zero.
\item The group $Q$ has the form $ANR$ where $A$ is simply connected and abelian,
$N$ is simply connected and unipotent, and $R$ is compact and containing the stabilizer
of $\omega_j$, $0\le j\le k$. In particular the stabilizer of $\omega_j$ in $Q$ is compact.
\end{enumerate}

Our aim in this article is to use the
special structure of the group $Q=ANR$ listed abouve, to construct frames
for $L^2(\mathbb{R}^n)$ generalizing (\ref{frame}). Our ideas are
based on the article \cite{BT96} by Bernier and Taylor , but it should be pointed out
that seveal of the ideas in \cite{BT96} are based on privious work
of others. We would like to mention the article
by  Daubechies, Grossmann, and  Meyer \cite{DGM86}, the
work of Feichtinger and Gr\"ochenig \cite{FG89,KG98}, Al, Antoine, and  Gazeau \cite{AAG2000,BT96},
and finally the work of Heil and Walnut \cite{HW89}.

In \cite{BT96} the authors considered a subgroup  $H\subset \GLR$ as above and
assume that $H$ acts
freely on $\R^n$, i.e., the stabilizer group is trivial. Define an action
of $H$ on $\mathbb{R}^n$ by $a\cdot x=(a^{-1})^T(x)$. The authors
introduced the notion  of \textit{separated} sets and \textit{frame generators}.
A separated set $\Gamma$ is a subset of $H$ such that there exist a compact set $B\subset \cO$,
$B^o\not= \emptyset$, such that $a\cdot B\cap b\cdot B\not=\emptyset$ implies that $a=b$.
In particular
all the translates $a^TB$, $a\in \Gamma$, of $B$ are disjoint. A \textit{frame generator} is a pair
$(\Gamma, \mathbb{F})$ where $\Gamma\subset H$ is separated and $\mathbb{F}$ is a compact
subset of $\cO$ such that $\bigcup_{a\in \Gamma} a\cdot \mathbb{F}=\cO$.
In \cite{BT96} the author show, that if $(\Gamma, \mathbb{F})$ is a frame generator, then
there exist a function $\psi$ and a discrete set $\{v(m)\in\R^n\mid m\in \Z^n\}$
such that the set
$\{\pi((a,v(m))^{-1})\psi\mid a\in \Gamma, m\in \Z^n\}$ is a
frame for the Hilbert space of $L^2$-functions whose Fourier transform
is supported in $\cO$.

In the first part of this article we show that the same construction can be
carried out if the action is not free, but the stabilizer group is compact.
Motivated by the construction in \cite{FO02} we apply this to groups
of the form $ANR$ as in (2) above except we do not need to assume that $N$ is
simply connected. The assumption needed is, that $N$ contains a discrete
subgroup $\Gamma$ such that $\Gamma \backslash N$ is compact. In this case
we can carry out the construction by Bernier and Taylor to get a frame related
to an open $ANR$-orbit, c.f. Theorem \ref{thmain1} and
Theorem \ref{thmain2}. We recall in section 4 the construction
from \cite{FO02} and explain it using the action of
$\GLR$ on the space of symmetric matrices.

Then author would like to thank C. Heil and G. Weiss for helpfull comments and
corrections.

\section{Separated sets}

Let $H$ be a closed subgroup of $\mathrm{GL}(n,\mathbb{R})$. Then $H$ acts in a
natural way on $\mathbb{R}^n$. We will also consider the action
$$(h,v)\mapsto h\cdot v:= (h^{-1})^T(v)$$
of $H$ on
$\mathbb{R}^n$. Here $a^T$ denotes the transpose of the matrix $a\in \mathrm{GL}(n,\mathbb{R})$.
We denote by $\theta :\R^n\to \R^n$ the homomorphism
$\theta (h)=(h^{-1})^T$. For simplicity we will also
write $h^\theta=\theta(h)$.
We assume that there
exist an open orbit $\cO\subset \R^n$ under the twisted action $(h,v)\mapsto
h^\theta(v)$. For $\omega\in\cO$ let
\begin{equation*}
H^{\omega }:=\{ h\in H\mid h \cdot \omega =\omega \}
\end{equation*}
be the \textit{stabilizer} of $\omega$ in $H$.
Notice that $H^{\omega }=\{h\in H\mid h^T(v)=v\}$ as $H^{\omega}$ is a subgroup
of $H$.
Because of the applications that we have in mind, we assume from now on
that $H^{\omega }$ is
compact.  The following definition
is from \cite{BT96}:

\begin{definition} Let $H$ be a locally compact Hausdorff topological group
acting on the locally compact Hausdorff topological space $X$.
A  subset $\Gamma \subset H$ is called \textbf{separated} if there exist a compact
set $B\subset X$ such that $B^{o}\not=\emptyset $ and $h\cdot B
\cap k\cdot B=\emptyset $ for all $h,k\in \Gamma $, $h\not=k$. We then say
that $\Gamma $ is separated by $B$.
\end{definition}

\begin{example}\label{SO1} Let $H=\R^+ \mathrm{SO}(n)$. Let $A\subset \mathrm{SO}(n)$ be
a non-empty subset and let $\lambda >1$.
$$\Gamma :=\{\lambda^k a\mid k\in \mathbb{Z},\, a\in A\}\, .$$
Let $0<\alpha <\beta$ be such that $\lambda\alpha > \beta$ and define
$$B=\{v\in \mathbb{R}^n\mid \alpha \le \|v\|\le \beta\}\, .$$
Then $B$ is compact with non-empty interior.
If $b=\lambda^k a\in \Gamma$ then
$$b\cdot B=\{v\in \mathbb{R}^n\mid \lambda^{-k} \alpha \le \|v\|\le \lambda^{-k} \beta \}\, .$$
Suppose that $k\le m$ and that $\lambda^k\mathrm{SO}(n)\cdot B\cap\lambda^m\mathrm{SO}(n)
\cdot B $.
Then $\lambda^{-k}\alpha \le \lambda^{-m }\beta $ and, hence,
$$\lambda^{m-k}\alpha \le  \beta $$
which is only possible if $m-k=0$. It follows that $\Gamma$ is separated by
$B$.
\end{example}

Fix from now on $\omega_0\in \mathcal{O}$ and recall that we are assuming
that $L:=\{h\in H\mid h^T( \omega_0)=\omega_0\}$ is compact.
We can always assume that $\omega_0 \in B^{o}$. Otherwise take $b\in H$ such that
$b\cdot \omega_0 \in B^{o}$. Thus $\omega_0=b^{-1}\cdot B$. Let $\Gamma^{\prime }:=\Gamma b $
and $B^{\prime }:=b^{-1}\cdot B$. Then
for $h,k\in \Gamma $, $h\not=k$ we have
\begin{equation*}
(hb)\cdot(b^{-1}\cdot B)\cap (kb)\cdot (b^{-1}\cdot B)
=
h\cdot B\cap k\cdot B
=\emptyset
\end{equation*}
so that $\Gamma^\prime$ is separated by $B^\prime$.

\begin{lemma} Let $L=\{h\in H\mid h^\theta (\omega_0)=\omega_0\}$ and let
$\kappa : H\to \mathcal{O}$ be the  map
$h\mapsto h^{\theta}( \omega_0)$ that defines an $H$-isomorphism $H/L\simeq \mathcal{O}$.
Let $B\subset \mathcal{O}$ be compact. Then $\tilde{B}:= \kappa^{-1}(B)\subset H$
is a right $L$-invariant compact subset of $H$ such that
$\kappa (\tilde{B})=B$.  Furthermore the
following holds:
\begin{enumerate}
\item $\omega_0\in B$ if and only if $e\in \tilde{B}$;
\item $B^o\not= \emptyset$ if and only if  $\tilde{B}^o\not= \emptyset$;
\item If $\tilde{B}^o\not=\emptyset$ then
$\kappa(\tilde{B}^o)=B^o$ and $\tilde{B}^o$ is right $L$-invariant.
\end{enumerate}
\end{lemma}
\begin{proof} All of this is well know, but let us
go over the argument here. That $\tilde{B}L=\tilde{B}$ follows from
the fact that $\kappa(ab)=a^\theta(\kappa(b))=(ab)\cdot \omega_0$.

(a) follows by $\kappa(e)=\omega_0$.

(b) and (c) follows from the fact that $\kappa$ is
open and continuous.

We can assume that $e\in \tilde{B}$. Then $L\subset \tilde{B}$.
Let $V\subset H$ be an open neighborhood of $e$ such
that $\bar{V}$ is compact. Then
$\kappa(V)\subset \mathcal{O}$ is open and
\[B\subset \bigcup_{g\in H}g\cdot \kappa(V)\, .\]
Hence, there are finitely many $g_1,\ldots ,g_n$ such that
\[B\subset \bigcup_{j=1}^n g_j\cdot \kappa(V)\, . \]
Let $S:=\cup_{j=1}^n g_j \overline{V}$. Then $S$ is compact and $B\subset \kappa(S)$. It
follows that $\tilde{B}\subset SL$. But $SL$ is compact as the continuous
image of the compact set  $S\times L\subset H\times H$
under the continuous map
$H\times H\to H$, $(a,b)\mapsto ab$. As $\tilde{B}$ is closed it follows that $\tilde{B}$ is
compact.
\end{proof}

\begin{lemma}\label{lemma1} Let $\Gamma\subset H$ be a separated set. Let $D\subset \mathcal{O}$ be compact.
Let $S=\kappa^{-1}(D)$. Then for
each $a\in \Gamma$:
\[ \#\{b\in \Gamma\mid a\cdot D\cap b\cdot D\not=\emptyset\}=\#\{b\in \Gamma \mid
aS\cap bS\not= \emptyset\}\, .\]
\end{lemma}
\begin{proof}
We have
\[\kappa(a S\cap bS)=a^\theta (D)\cap b^\theta (D)\, .\]
Hence, if $aS\cap bS\not=\emptyset$ then $a^\theta (D)\cap b^\theta (D)\not=\emptyset$ and it
follows that the right hand side is greater or equal to the left hand side. Assume now
that $x\in a^\theta (D)\cap b^\theta (D)$. Then there exists $s,t\in S$ such that
$\kappa (as)=\kappa (bt)$. Hence, there exist $h\in L$ such that
$as=bth$. As $SL=S$ it follows that $aS\cap bS\not=\emptyset$ and,
hence,
$$\#\{b\in \Gamma\mid a\cdot D\cap b\cdot D\not=\emptyset\}\le \#\{b\in \Gamma \mid
aS\cap bS\not= \emptyset\}$$
finishing the proof.
\end{proof}
We have now the necessary tools to prove the main results of this section.

\begin{theorem}\label{lemma2}
Suppose that $\Gamma \subset H$ and that $B\subset \cO$. Let $\tilde{B}=\kappa^{-1}(B)$.
Then $\Gamma$ is separated by $B$ in $\cO$ if and only of
$\Gamma$ is separated by $\tilde{B}$ in $H$  under the natural action of
$H$ on $H$ given by left multiplication.
\end{theorem}

\begin{proof} We have already seen that $\tilde{B}$ is compact with $\tilde{B}^o$ non-empty.
Assume that $a\tilde{B}\cap b\tilde{B}\not=\emptyset$ for some $a,b\in \Gamma$
then it follows by Lemma \ref{lemma1} that $a\cdot B\cap b\cdot B\not=\emptyset$.
Hence, $a=b$ as $\Gamma$ is separated by $B$.
\end{proof}

\begin{theorem}\label{leFinite}
Suppose that $\Gamma \subset H$ is a separated subset of $H$. Let $D\subset \mathcal{%
O}$ be compact. Then
\begin{equation*}
\sup_{k\in \Gamma }\#\{h\in \Gamma \mid h \cdot D\cap k\cdot D\not=\emptyset
\}<\infty \,.
\end{equation*}
\end{theorem}

\begin{proof}
This has been proved in \cite{BT96} for the case where
the action of $H$ is
free. Using Lemma \ref{lemma1} and Theorem \ref{lemma2}
the general statement is reduced to that case and, hence, the claim.
\end{proof}

\section{The continuous wavelet transform}
In this section we review some basic facts about the continuous wavelet transform
on $\mathbb{R}^n$ with respect to a group action, 
see \cite{FO02,GMP85,GMP86} for more information and references. Denote by $\mathrm{Aff}(\R^n)$ the group of
invertable affine linear
transformations on $\mathbb{R}^n$. Then $\mathrm{Aff}(\R^n)$ consists of pairs $(x,h)$ such
that $h\in \GLR$ and $x\in \mathbb{R}^n$. The action of
$(x,h)\in \mathrm{Aff}(\R^n)$ on $\mathbb{R}^n$ is given by
\[(x,h)(v)=h(v)+x\, .\]
The product is the composition of maps. Thus
\[(x,a)(y,b)=(a(y)+x,ab)\]
and the inverse of $(x,a)\in \mathrm{Aff}(\R^n)$ is given by
\[(x,a)^{-1}=(-a^{-1}(x),a^{-1})\, .\]
Thus $\mathrm{Aff}(\R^n)$ is the semidirect product of the abelian group $\mathbb{R}^n$ and
the group $\mathrm{GL}(n,\mathbb{R})$;
$\mathrm{Aff}(\R^n)=\mathbb{R}^n\times_s \mathrm{GL}(n,\mathbb{R})$.

Define a unitary representation of $\mathrm{Aff}(\R^n)$ on $L^2(\mathbb{R}^n)$ by
\begin{equation}
[\pi (x,a)f] (v)=|\det(a)|^{-1/2}f((x,a)^{-1}(v)
=|\det(a)|^{-1/2}f(a^{-1}(v-x))\, .
\end{equation}
For $f\in L^2(\mathbb{R}^n)$ denote
by $\hat{f}$ the Fourier transform of $f$
\[\hat{f}(\omega) =\frac{1}{(2\pi)^{n/2}}\int_{\mathbb{R}^n} f(x)e^{ - i(x\mid \omega)}
\, dx\, ,\quad f\in L^1(\R^n)\cap L^2(R^n)\, .\]
We denote by $\hat{\pi}(x,a)$ the unitary action on $L^2(\mathbb{R}^n)$ given by
\[\hat{\pi}(x,a)f(v)=\sqrt{|\det (a)| }e^{-i(x\mid v)}f(a^T(v))=\sqrt{|\det (a)|}e^{-i(x\mid v)}f(a^{-1}\cdot v)\, .\]
The Fourier transform intertwines the representations $\pi$ and $\hat{\pi}$ \cite{FO02}, Lemma 3.1:
\begin{lemma}\label{le-actionFour} Let $f\in L^2(\mathbb{R}^n)$ and $(x,a)\in \mathrm{Aff}(\R^n)$.
Then
\[\widehat{\pi (x,a)f}(\omega)= \hat{\pi}(x,a)\hat{f}(\omega )\, .\]
\end{lemma}

Let $H\subset \mathrm{GL}(n,\mathbb{R})$ be a closed subgroup. Denote by
$G:=\mathbb{R}^n\times_s H$ the subgroup of $\mathrm{Aff}(\R^n)$ given by
\[G=\{(x,a)\in \mathrm{Aff}(\R^n)\mid a\in H\}\, .\]
We assume that there exists open  sets $\{\mathcal{O}_j\}_{j\in \mathbb{J}}$, where
$\mathbb{J}$ is a finite or countably infinite index set, such that
\begin{enumerate}
\item[(W1)] Each $\mathcal{O}_j$ is invariant and homogeneous under the action of $H$ given by
$(a,v)\mapsto a\cdot v=a^\theta(v)$;
\item[(W2)] We have $\mathcal{O}_i\cap \mathcal{O}_j=\emptyset$ if $i\not= j$;
\item[(W3)] The complement of $\cup_{j\in \mathbb{J}}\mathcal{O}_j$ has measure zero
with respect to the Lebesgue measure on $\mathbb{R}^n$.
\end{enumerate}
For a measureable function $f$ denote by $\mathrm{Supp}(f)$ the complement of the maximal open set $U\subset \mathbb{R}^n$
such that $f(x)=0$ for almost all $x\in U$.  For $\emptyset\not= U$ open in $\mathbb{R}^n$ denote
by $L^2_U(\mathbb{R}^n)$ the closed subspace of $L^2(\mathbb{R}^n)$ given by
\begin{equation}\label{L2U}
L^2_U(\mathbb{R}^n)=\{f\in L^2(\mathbb{R}^n)\mid
\mathrm{Supp}(\hat{f})\subset \overline{U}\}\, .
\end{equation}
Suppose that $U$ is $H$-invariant under the twisted action $h\cdot v=h^\theta(v)$, then by \ref{le-actionFour}, it follows
that $L^2_U(\mathbb{R}^n)$ is $G$-invariant. Furthermore, by Theorem 3.4 in \cite{FO02},
$L^2_U(\mathbb{R}^n)$ is
irreducible if and only if $U$ is homogeneous.
Hence, the decomposition of $L^2(\mathbb{R}^n)$ into irreducible parts is given by
\begin{equation}\label{irred}
L^2(\mathbb{R}^n)\simeq \bigoplus_{j\in \mathbb{J}}L^2_{\mathcal{O}_j}(\mathbb{R}^n)\, .
\end{equation}
Denote by $dh$ a left invariant Haar measure on $H$. Then a
left invariant Haar measure on $G$ is given by $dg=(2\pi )^{-n}|\det (a)|^{-1}dadv$.

\begin{definition} Suppose that $\emptyset\not= U$ is
an open subset of $\mathbb{R}^n$. Then a nonzero function
$f\in L^2_{\mathcal{O}_j}(\mathbb{R}^n)$ is called
\textbf{admissible}  if
$\pi_{g,f}(h):=(g\mid \pi (x,a)f)$ is in $L^2(G)$ for
all $g\in L^2_{\mathcal{O}_j}(\mathbb{R}^n)$.
\end{definition}

A simple calculation, see \cite{FO02}, shows that
\begin{equation}\label{wavelet}
\int_G|(g\mid \pi (x,a)f)|^2\frac{dadx}{|\det (a)|}=
(2\pi )^n \int_U |\hat{g}(\omega )|^2\int_{H}|\hat{f}(h^T(\omega ))|^2\, dhd\omega
\end{equation}
In particular, if $U$ is homogeneous, then $C_f=\int_{H}|\hat{f}(h^T(\omega ))|^2\, dh$ is
independent of $\omega\in U$ and, hence,
\[\int_G|(g\mid \pi (x,a)f)|^2\frac{dadx}{|\det (a)|}=C_f\|g\|^2\, .\]
In particular $f$ is admissible if and only if $H\ni h \mapsto f(h^T\omega)\in\C$ is in $L^2(H)$, which in
particular implies that the condition
\begin{enumerate}
\item[(W4)] For all $\omega \in U$ we have that
$H^\omega =\{h\in H\mid h^T(\omega)=\omega\}$ is compact
\end{enumerate}
has to be satisfied.

\section{Separated sets and Frames}

In this section we recall some basic facts from
\cite{BT96} on how to  construct frames from the continuous
wavelet transform using separating sets. Let us also recall
that we are assuming that $H\subset \mathrm{GL}(n,\mathbb{R})$ is
closed and that $\mathcal{O}$ is a homogeneous open subset of $\mathbb{R}^n$
such that the condition (W4) is satisfied.

Let us start with the well known definition:

\begin{definition}
Let $\mathbf{H}$ be a
Hilbert space. A sequence $\{v_n\}$ in $\mathbf{H}$ is called a \textbf{frame}
if there exits numbers $A,B>0$ such that for all $v\in \mathbf{H}$ we have
\[A\|v\|^2\le \sum_n |(v\mid v_n)|^2\le B\|v\|^2\, .\]
The numbers $A$ and $B$ are called \textbf{frame bounds}.
\end{definition}

The following definition is a simple generalisation of the definition by
Bernier and Taylor \cite{BT96}:

\begin{definition} Let $H$ be a locally compact Hausdorff topological group
acting transitively on the locally compact Hausdorff topological space $X$.
A \textbf{frame generator} is a pair $(\Gamma , \mathbb{F})$ where $\Gamma$
is a countable  separated subset of $H$ and $\mathbb{F}$ is a compact subset of $X$ such
that
\begin{equation}
X=\bigcup_{a\in \Gamma} a\cdot \mathbb{F}
\, .\label{framegen}
\end{equation}
\end{definition}

In our case we will take $X=\cO$ or $X=H$. Notice (\ref{framegen}) implies in this case
that $\mathbb{F}^o\not=\emptyset$ and that for each $a\in \Gamma$
we have $\#\{b\in \Gamma\mid a\cdot \mathbb{F}\cap b\cdot \mathbb{F} \not=\emptyset\}<\infty$.
Notice also that $(\Gamma, \mathbb{F})$ is a frame generator for the action on $\cO$
if and only if $(\Gamma ,\kappa^{-1}(\mathbb{F}))$ is a frame generator for the
action of $H$ on $H$ by multiplication.

Let us now go back to the situation considered in the
previous sections. Let $(\Gamma ,\mathbb{F})$ be a frame generator, let $D\subset \mathcal{O}$ be
compact subset of $\mathcal{O}$ such that $\mathbb{F}\subset D^o$. Let $R\subset \mathbb{R}^n$
be a parallelepiped such that $D\subset R$. Choose $a_j<b_j$ ($j=1,\ldots ,n$) and
a basis $v_j\in \mathbb{R}^n$ ($j=1,\ldots ,n$) such that
\[R=\{\sum_{j=1}^n x_jv_j\mid a_j\le x_j\le b_j\}\, .\]
Let $w_1,\ldots ,w_n$ be
the dual base to $v_1,\ldots ,v_n$, i.e., $(v_i\mid w_j)=\delta_{ij}$. For
$m=(m_1,\ldots ,m_n)\in \mathbb{Z}^n$ define $w(m)\in \mathbb{R}^n$ by
\[w(m):=\sum\frac{m_j}{b_j-a_j}\, w_j\, .\]
Finally we define $e_m:\mathbb{R}^n\to \mathbb{C}$ by
\[e_m(v) =\frac{1}{\sqrt{\mathrm{Vol}(R)}}
 \exp(\sum_{j=1}^n 2\pi i (v\mid w(m)))\chi_{R}(v)\]
 where $\chi_F$ denotes the indicator function of a set $F\subset \mathbb{R}^n$.
 We identify $e_m$ with its restriction to $R$.
Then $\{e_m\}_{m\in \mathbb{Z}^n}$ is an orthonormal basis for $L^2(R)$. Let
\[\alpha :=\sup_{a \in \Gamma} \#\{b\in \Gamma\mid a\cdot  D\cap b\cdot D\not=
\emptyset\}\, . \]
Then $\alpha$ is finite by Lemma \ref{leFinite}.
Let $\varphi \in L^2_{\mathcal{O}}(\mathbb{R}^n)$
be such that:
\begin{enumerate}
\item[(F1)] $\mathrm{Supp}(\hat{\varphi})\subset D$;
\item[(F2)] $a(\varphi ):=\inf_{\omega \in \mathbb{F}}|\hat{\varphi}(\omega )|>0$;
\item[(F3)] $b(\varphi):=\sup_{\omega \in D}|\hat{\varphi}(\omega )|<\infty$.
\end{enumerate}
In particular we could take $\varphi$ such that $\hat{\varphi}=\chi_{\mathbb{F}}$.

\begin{theorem}[Bernier-Taylor] Assume that $(\Gamma, \mathbb{F})$ is
a frame generator and that $\varphi\in L^2_{\mathcal{O}}(\mathbb{R}^n)$
satisfies the conditions (F1), (F2), and (F3). Then, with the above notation, the sequence
\[\{\pi( (a ,w(m))^{-1})\varphi \}_{(a,m)\in \Gamma\times \mathbb{Z}^n}\]
is a frame for $L^2_{\mathcal{O}}(\mathbb{R}^n)$ with frame bounds
$A=\mathrm{Vol}(R)a(\varphi )^2$ and $B=
\mathrm{Vol}(R)\alpha b(\varphi)^2$.
\end{theorem}

\begin{proof} See \cite{BT96}, Theorem 3.
\end{proof}

\begin{example}\label{SO2} Let $H=\R^+ \mathrm{SO}(n)$ as in Example \ref{SO1}.
Then $H$ has two orbits, $\{0\}$ and $\mathcal{O}:=\mathbb{R}^n\setminus \{0\}$
in $\mathbb{R}^n$.  Notice that
$L^2_{\mathcal{O}}(\mathbb{R}^n)=L^2(\mathbb{R}^n)$ in this
case. If $u\in \mathcal{O}$ and $r\not= 1$ then $ru\not= u$ and,
hence, $H^\omega$ is a closed subgroup of $\mathrm{SO}(n)$ and therefore compact.
In fact it is easy to see that $H^\omega$ is isomorphic to $\mathrm{SO}(n-1)$
for all $\omega\in \mathcal{O}$. Thus all the conditions (W1) -- (W4) are fullfilled.

Let
$\lambda >1$ and let $\Gamma =\{\lambda^n\mid n\in \mathbb{Z}\}$. Then $\Gamma$
is separated by example \ref{SO1}. Choose $\rho <\sigma$ such that $\lambda \rho \le \sigma $ and define
$\mathbb{F} = \{v\in \mathbb{R}^n\mid \rho \le \|v\|\le \sigma\}$. Let
$v\in \cO$. Then there exist an $n\in \mathbb{Z}$ such
that $\lambda^{n}\rho \le \|v\|<\lambda^{n+1}\rho\le \lambda^{n}\sigma$.
Hence, $v\in \lambda^{-n}\mathbb{F}=\Gamma\cdot \mathbb{F}$. Thus $(\Gamma,\mathbb{F})$ is a
frame generator.
\end{example}

\begin{example}\label{HinSO1n} In this example we consider a case of
a subgroup $H\subset \R^+ \mathrm{SO}(1,n)$ acting on $\mathbb{R}^{n+1}$
which is more complicated than the example \ref{SO2}. But our
construction relies on the fact that the action on
each of the open orbits is free, i.e., the stabilizer
is trivial.

For $\lambda\in\R^*=\{r\in\R\mid r\not= 0\}$, $t\in \mathbb{R}$,
and $x\in \R^{n-1}$ define $a(\lambda,t),n(x)\in \GL (n+1,\R)$ by
\[a(\lambda ,t)=\lambda \left(
\begin{matrix}
\cosh(t) & 0 & \sinh(t)\\
0 & \mathrm{I}_{n-1} & 0\\
\sinh(t) & 0 & \cosh(t)
\end{matrix}\right)\]
and
\[
n(x)=\left(
\begin{matrix}
1+\oh \|x\|^2& x^T & \frac{1}
{2}\|  x \|^{2}\\
x & \mathrm{I}_{n-1} & x\\
- \frac{1}{2}\|  x\|^{2} & -x^T & 1-\frac{1}
{2}\|  x \|  ^{2}
\end{matrix}
\right)\]
and let
\[A:=\{a(\lambda,t)\mid \lambda >0, t\in\mathbb{R}\}
\quad \mathrm{and}
\quad N=\{n(x)\mid x \in \mathbb{R}^{n-1}\} \, .\]
Then $A$ and $N$ are abelian groups. But calculations
are in fact easier using the corresponding
Lie algebras, which are abelian and isomorphic to
$\R^2$, respectively $\R^{n-1}$. For that let
\begin{equation*}
H(s,t)=s I_{n+1}+s(E_{1\, n+1}+E_{n+1\, 1})\quad \mathrm{and}\quad
X(x)=\left(\begin{matrix} 0 & x^T & 0\cr
x & 0 & x\cr
0 & -x^T & 0\end{matrix}\right)\, ,
\end{equation*}
where $x\in \mathbb{R}^{n-1}$ and $E_{\nu \mu }=(\delta_{i\nu }\delta_{j\mu})_{ij}$.
Define
\begin{equation*}
\mathfrak{a}=\{H(s,t)\mid s,t\in \mathbb{R}\}\simeq \mathbb{R}^2\quad \mathrm{and}
\quad
\mathfrak{n}=\{X(x)\mid x\in \mathbb{R}^{n-1}\}\simeq \mathbb{R}^{n-1}\, .
\end{equation*}
Denote by
$$X\mapsto \exp (X)=e^X=\sum_{j=0}^\infty \frac{X^j}{j!}$$
the matrix exponential function. Then
$$e^{H(s,t)}=a(e^s,t)\quad \mathrm{and}\quad
e^{X(x)}=n(x)\, .$$
Furthermore $\exp : \mathfrak{a}\to A$ and $\exp :\mathfrak{n}\to N$ is a
group homomorphism, i.e., in both cases we have $e^{X+Y}=e^Xe^Y$ ($X,Y\in \fa$ or $X,Y\in\fn$),
hence, the multiplication in $A$, respectively, $N$ can be reduced to
the usual addition in $\R^2$, respectively $\R^{n-1}$.

A simple calculation shows that
\begin{equation*}
a(\lambda ,t)n(x)a (\lambda ,t)^{-1}=n(e^{-t}x )\, .
\end{equation*}
In particular it follows that
$H=AN=NA$
is a closed subgroup of $\mathrm{GL}(n+1,\mathbb{R})$
with $N$ a normal subgroup. Next we notice that
$a(\lambda ,t)^\theta =a(\lambda ,t)^{-1}= a(\lambda^{-1} , -t)$, $n(x)^{-1}=n(-x)$, and
\[n(x)^{\theta}=n(-x)^{T}=\left(\begin{matrix}
1+\frac{1}{2}\|x\|^2 & -x^T & -\frac{1}{2}\|x\|^2\cr
-x & I_{n-1} & x\cr
\frac{1}{2}\|x\|^2 & -x^T & 1-\frac{1}{2}\|x\|^2
\end{matrix}\right)\, .\]
Hence, the twisted action of $a(\lambda, t)$ and $n(x)$ is
given by
$$a(\lambda ,t)\cdot v=\lambda^{-1} (\cosh (t)v_1-\sinh(t)v_{n+1},v_2,\ldots ,v_{n},
-\sinh(t)v_1+\cosh (t)v_{n+1})^T$$
and
$$n(x)\cdot v=v+(v_1- v_{n+1})(\frac{1}{2}\|x\|^2, -x^T, \frac{1}{2}\|x\|^2)^T-
\left(\sum_{j=1}^{n-1}x_jv_{j+1}\right)(1,0,\ldots ,0,1)^T\, .$$
In particular, if we take $v= e_1$ and $v= e_{n+1}$, where $\{e_j\}_{j=1}^{n+1}$ is
the standard basis of $\mathbb{R}^n$, we get
\[n(x)a(\lambda ,t)\cdot e_1=
\lambda^{-1} (\cosh(t)+\frac{e^t}{2}\|x\|^2,-e^t x^T, -\sinh(t)+\frac{e^t}{2}\|x\|^2)^T\]
and
\[n(x)a(\lambda ,t)\cdot e_{n+1}=\lambda^{-1}(-\sinh (t)-\frac{e^{t}}{2}\|x\|^2,e^t x^T,\cosh (t)-\frac{e^{-t}}{2}\|x\|^2)^T\, .\]
Notice that the stabilizer of $e_1$ and $e_{n+1}$ is trivial. Let $\beta$ be the bilinear form
\[\beta (v,w)=v_1w_1-\sum_{j=2}^{n+1}v_jw_j\, .\]
We can now describe the four open $H$-orbits:

\begin{eqnarray*}
\mathcal{O}_1& =&\{v\in \mathbb{R}^{n+1}\mid \beta (v,v)>0, v_1>0\}=H\cdot e_1\\
\mathcal{O}_2& =&\{v\in \mathbb{R}^{n+1}\mid \beta (v,v)>0, v_1<0\}=H\cdot (-e_1)\\
\mathcal{O}_3& =&\{v\in \mathbb{R}^{n+1}\mid \beta (v,v) < 0, v_1<v_{n+1}\}=H\cdot e_{n+1}\\
\mathcal{O}_4& =&\{v\in \mathbb{R}^{n+1}\mid \beta (v,v)<0, v_1>v_{n+1}\}=H\cdot (-e_{n+1})
\end{eqnarray*}
The complement of $\mathcal{O}_1\cup \cdots \cup \mathcal{O}_4$ is given by
\[\mathbb{R}^{n+1}\setminus (\mathcal{O}_1\cup \cdots \cup \mathcal{O}_4)
=\{v\in \mathbb{R}^{n+1}\mid \beta (v,v)=0\}
\cup \{v\in \mathbb{R}^{n+1}\mid \beta (v,v)<0,\, v_1=v_{n+1}\}
\]
which obviously has measure zero in $\mathbb{R}^n$. Thus (W1) -- (W4) holds. Notice also, that
if we allow $\lambda$ to take positive and negative values, i.e., replace $H$ by the non-connected group
$\{a(\lambda ,t)\mid \lambda\in \R^*,\, t\in \R\}N$; then there are only two open orbits,
$\cO_1\cup \cO_2$ and $\cO_3\cup \cO_4$.

Define now
$$\Gamma_A=\{\exp (H(n,m))\mid n,m\in \mathbb{Z}\}
\quad \mathrm{and}\quad \Gamma_N=\{\exp X(x)\mid x\in \mathbb{Z}^n\}\, .$$
Then $\Gamma_A$ and $\Gamma_N$ are discrete subgroups of $A$ respectively $N$ and $N/\Gamma_N$ is
compact. Let
\begin{equation}\label{defgamma}
\Gamma :=\Gamma_A\Gamma_N\subset H\, .
\end{equation}
Then $\Gamma$ is a discrete subset of $H$, but notice that $\Gamma$ is not a group.
For $\epsilon >0$ denote by
$$B_{\epsilon}=\{x\in \mathbb{R}^{n-1}\mid
\|x\|\le \epsilon\}\, .$$
Choose $0<\delta<1/4$ such that $e^tB_{1/2}\subset B_{3/4}$ for all $t\in [-\delta ,\delta]$.
Then $e^{t}\mathbb{Z}^{n-1}\cap B_{1/2}=\{0\}$ for all $|t|\le \delta$.
Let
$$B=\{a(e^s,t)n(x)\mid |s|\le 1/4,\, |t|\le \delta,\, x\in B_{1/2}\}\, .$$
Then $B\subset H$ is  compact and $B^o\not=\emptyset$. Let $a,b
\in\Gamma$ and assume
that $aB\cap bB\not= \emptyset$. Then there exists $a(e^r,t)n(x), a(e^s,u)n(y)\in B$, such
that $aa(e^r,t)n(x)=ba(e^s,u)n(y)$. Write $a=a(e^{n_1},m_1)n(\mathbf{m}_1)$
and $b=a(e^{n_2},m_2)n(\mathbf{m}_2)$ with $n_j,m_j\in \mathbb{Z}$, and $\mathbf{m}_j
\in \mathbb{Z}^{n-1}$ ($j=1,2$). Then
we have
\begin{eqnarray*}
aa(e^r,t)n(x)&=&a(e^{n_1+r},m_1+t)n(e^{r}\mathbf{m}_1+x)\\
&=& ba(e^s,u)n(y)\\
&=&a(e^{n_2+s},m_2+u)n(e^{s}\mathbf{m}_2+y)\, .
\end{eqnarray*}
But this is only possible if
\[n_1+r=n_2+s,\quad m_1+t=m_2+s\, \quad\mathrm{and}
\quad e^{r}\mathbf{m}_1+x=e^{s}\mathbf{m}_2+y\, .
\]
But then $n_1-n_2=s-r\in \mathbb{Z}\cap [-1/2,1/2]=\{0\}$ and, hence, $n_1=n_2$ and $s=r$. Similarly
it follows that $m_1=m_2$ and $t=u$. Thus
$$e^{r}\mathbf{m}_1+x=e^{r}\mathbf{m}_2+y$$
or
$$e^{r}(\mathbf{m}_1-\mathbf{m}_2)=y-x\in e^{r}\mathbb{Z}^{n-1}\cap B_{1/2}=\{0\}$$
which implies that $\mathbf{m}_1=\mathbf{m}_2$ and $y=x$. In particular it follows
that $\Gamma$ is separated by $B$ in $H$.
Let $\omega_1=e_1$, $\omega_2=-e_2$, $\omega_3=e_{n+1}$ and $\omega_4=-e_{n+1}$.
For $j=1,2,3,4$ the map
$$H\ni h\mapsto \kappa_j(g):=h^{\theta}(\omega_j)=h\cdot \omega_j\in \mathcal{O}_j$$
is a diffeomorphism such that $\kappa_j(ab)=a\cdot \kappa_j(b)$. It
follows that $\Gamma$ is separated by $B_j:=\kappa_j(B)$ ($j=1,\ldots ,4$).

For $j=1,\ldots ,r$ let
$$
\mathbb{F}_j=\kappa_j(\{a(e^s,t)n(x^T)\mid |s|\le 1, |t|\le 1, |x_k|\le 1 \, (k=1,\ldots ,n-1)\})\, .
$$
Then a simple calculation shows that $(\Gamma,\mathbb{F}_j)$ is a frame generator.
\end{example}

\begin{example}\label{ex3.6} Let $H$ be a locally compact Hausdorff topological group, and assume that there exist a
countable discrete subgroup
$\Gamma\subset H$ such that $\Gamma\backslash H$ is compact.
Then there exist a compact subset $K\subset H$ such that $e\in K^o$,
$K^{-1}=\{k^{-1}\mid k\in K\} =K$ and $K^2=\{ab\mid a,b\in K\}\cap \Gamma =\{e\}$.
Assume that $a,b\in \Gamma$ and $aK\cap bK\not=\emptyset$. Then
$b^{-1}a\in K^2\cap \Gamma=\{e\}$ and therefore $a=b$. It follows that $\Gamma$ is separated
by $K$. As $\Gamma\backslash H$ is compact there exist $\mathbb{F}\subset H$ such
that $\Gamma \mathbb{F}=H$. Thus $(\Gamma ,\mathbb{F})$ is a frame generator.

Assume now that $H\subset \GLR$ is closed and that $\cO\subset \R^n$ is an open
orbit. We assume that there exist $\omega\in \cO$ such that $\Gamma\cap H^\omega=\{e\}$.
Let, as usual, $\kappa :H\to \cO$ be the canonical map $\kappa (a)=a^\theta (\omega )
=(a^{-1})^T(\omega)$. Then $\Gamma$ is separated by $B=\kappa (K)$ and
$(\Gamma ,\kappa (\mathbb{F}))$ is a frame generator. We would expect that by modifying
the proof of Theorem 3 in \cite{BT96} one can remove the condition that
$\Gamma \cap H^\omega =\{e\}$, which would give several examples of reductive groups
acting on $\R^n$. The question is also, if one can remove the condition that
$H^\omega$ is compact, by $\Gamma \cap H^\omega$ is finite and
one assumes that $\Gamma \backslash G/H^\omega$ is compact.
\end{example}

\section{Action of some special groups and frame generators}
\label{sec3}
There are natural examples where
$\mathbb{R}^n$ contains finitely many open orbits satisfying (W1) -- (W3) but some of which do not
have compact stabilizers. In \cite{FO02} it was shown that in the case
where $H$ is reductive, or more simply stated,  $H^T=H$ and the stabilizer $L=H^\omega$, $\omega\in \cO$, is a symmetric
subgroup of $H$, then one can always find a subgroup $Q=RAN=ANR$ such that $\cO$ decomposes -- up to a set
of measure zero --
into finitely many open orbits such that (W1) -- (W4) holds. More importantly, the structure of
the group $Q$ is relatively simple and well understood. In particular we have the following:
\begin{enumerate}
\item The map $A\times N\times R\ni (a,n,r)\mapsto anr\in Q$ is a
diffeomorphism;
\item $R$ is a compact group and
the stabilizer of $\omega$ is contained in $R$;
\item $A$ is abelian and $A$ and $R$ commutes;
\item Both $R$ and $A$ normalize the group $N$;
\item Let
\[\fa=\{X\in M(n,\R)\mid \forall t\in\R\, :\, e^{tX}\in A\}
\]
be the Lie algebra of $A$.
Then $\exp : \fa\to A$ is an isomorphism of abelian groups, i.e., $e^{X+Y}=e^Xe^Y$ for $X,Y\in\fa$;
\item Let
\[\fn=\{X\in M(n,\R)\mid \forall t\in\R\, :\, e^{tX}\in N\}\]
be the Lie algebra of $N$. Then $\exp :\fn \to N$ is a diffeomorphism.
\end{enumerate}

We refer to \cite{FO02} for the exact construction, but for completeness we
recall some of the main constructions in the next section, but first we will show how those
facts can be used to construct a separated set and a frame generator.
The construction is very much in the spirit of Example \ref{HinSO1n}.
In fact we do not need the last conditions, so from now on we will assume that $Q=ANR\subset \GLR$
is a closed subgroup satisfying the
conditions (1) -- (5).

Let $H_1,\ldots ,H_r$ be a basis for $\fa$. For $\mathbf{t}\in \R^r$
let $a(\mathbf{t}):=\exp (\sum_{j=1}^rt_jH_j)$. Let
\begin{equation}\label{GammaA}
\Gamma_A=\{a(\mathbf{m})\in A\mid \mathbf{m}\in\mathbb{Z}^r\}\, .
\end{equation}
Then $\Gamma_A$ is a discrete subgroup of $A$. If $\Gamma_N$ is
a discrete subgroup of $N$ let $\Gamma=\Gamma_A\Gamma_N$. Then
$\Gamma$ is a discrete subset of $AN$, but in general $\Gamma$ is
not a subgroup.

\begin{theorem}\label{thmain1} Assume that $Q=ANR$ is a closed subgroup
of $\GLR$ such that (1) -- (5) above holds.
Suppose that $\Gamma_N$ is a discrete subgroup of $N$. Set $\Gamma=\Gamma_A
\Gamma_N$.   Assume that $\cO\subset \R^n$ is
an open  $Q$-orbit such that
$\overline{\cO}\setminus \cO$ has measure zero. Let $\omega\in \cO$ and assume
that $Q^\omega=L \subset R$. Then $\Gamma$ is a separated set.
\end{theorem}

\begin{proof} For $\epsilon >0$ denote by $B_\epsilon$ the closed ball in $\fn$
with center zero and radius $\epsilon$, $B_\epsilon =\{X\in \fn\mid
\mathrm{Tr}(XX^T)\le \epsilon\}$. Let $K_\epsilon =\exp (B_\epsilon)$. Choose $\epsilon >0$ such that
$\exp : B_{2\epsilon}^o\to \exp (B_{2\epsilon}^o)$ is a diffeomorphism, and
\begin{equation}\label{Keps1}
K_\epsilon^4 \cap \Gamma=\{e\}.
\end{equation} Choose $0<\delta\le 1/4$
such that for $|t_j|\le \delta$, $j=1,\ldots, r$, we have
\begin{equation}\label{Keps2}
\exp (a(\mathbf{t}) B_{\epsilon}a(\mathbf{t})^{-1})=a(\mathbf{t})K_\epsilon a(\mathbf{t})^{-1} \subset K_\epsilon^2
\, .
\end{equation}
This is possible because the action $A\times \fn \ni (a,X)\mapsto aXa^{-1}\in\fn $ is continuous.
For $X\in\fn$ let $n(X)=\exp (X)$.
Define (using the obvious notation)
\[B(Q)=\{a(\mathbf{r})n(X)b\mid |r_j|\le \delta ,\, X\in B_\epsilon,\, b\in B\}\subset Q\, .\]
Then $B(Q)$ is compact with $B(Q)^o\not=\emptyset $.

Assume that we have $g_1=\gamma_1\eta_1,g_2=\gamma_2\eta_2\in \Gamma$, $\gamma_j\in \Gamma_A$ and $\eta_j\in \Gamma_N$ ($j=1,2$),
such that $g_1B(Q)\cap g_2B(Q)\not=\emptyset$.
Then we can find $a_j=a(\mathbf{r}_j)\in \{a(\mathbf{r})\mid |r_j|\le \delta \}\subset A$,
$n_j=n(X_j)\in \{n(X)\mid  X\in B_\epsilon\}$, and $b_j\in B$ ($j=1,2$), such that
\[\gamma_1\eta_1a_1n_1b_1=\gamma_2\eta_2a_2n_2b_2\, .\]
But then
\[\gamma_1a_1((a_1^{-1}\eta_1a_1)n_1)b_1=\gamma_2a_1(a_2^{-1}\eta_2a_2)n_2b_2\, .\]
As the map $Q\simeq A\times N\times R$ (cf. condition (a)) we must have
$\gamma_1a_1=\gamma_2a_2$, $(a_1^{-1}\eta_1a_1)n_1=(a_2^{-1}\eta_2a_2)n_2$, and
$b_1=b_2$. Write $\gamma_j=a(\mathbf{m}_j)$ with $\mathbf{m}_j\in \mathbb{Z}^r$
 ($j=1,2$). Then
\[ a(\mathbf{m}_1+\mathbf{r}_1)=
a(\mathbf{m}_2+\mathbf{r}_2)\, .\]
As the exponential map $\exp :\fa\to A$ is an isomorphism of groups (cf. (5)) it follows that
$\mathbf{n}_1-\mathbf{n}_2=\mathbf{r}_2-\mathbf{r}_1\in \mathbb{Z}^r\cap \{\mathbf{r}\in
\mathbb{Z}^r\mid |r_j|\le \delta \}=\{0\}$.
Hence, $\mathbf{n}_1=\mathbf{n}_2$ and $\mathbf{r}_1=
\mathbf{r}_2$.
It follows that $\gamma_1=\gamma_2$ and
$a_1=a_2$. Let $a=a_1=a_2$. Then we have
\[a^{-1}\eta_1an(X_1)=a^{-1}\eta_2an(X_2)\]
or (by conjugating by $a$):
\[\eta_1 n(aX_1a^{-1})=\eta_2n(aX_2a^{-1})\, .\]
Write $\eta=\eta_2^{-1}\eta_1\in \Gamma_N$. Then we get
\[\eta=n(aX_2a^{-1})n(-aX_2a^{-1})\in \Gamma\cap K_\epsilon^4=\{e\}\]
by (\ref{Keps1}) and (\ref{Keps2}). But then
$\eta_1=\eta_2$ and $n_1=n_2$ showing that $B(Q)$ separate $\Gamma$ in $Q$.

Let $\kappa :Q\to \cO$, $p\mapsto p\cdot \omega=(p^{-1})^T(\omega)$.
Then $\kappa (pq)=p\cdot \kappa (q)$. Then $B:=\kappa (B(Q))$ is compact and
$B^o=\kappa (B(Q)^o)\not=\emptyset$. We claim that $B$ separate $\Gamma$ in $\cO$.
For that assume that there are $g_1,g_2\in \Gamma$ such that
$g_1\cdot B\cap g_2\cdot B\not=\emptyset$. Then there exists $b_1,b_2\in B(Q)$ such
that $g_1b_1=g_2b_2$. Then it follows that $g_1=g_2$ and, hence, the claim.
\end{proof}

Assume now that $\Gamma_N$ is a discrete subgroup of $N$ such that
$\Gamma_N\backslash N $ is compact. Choose $F_N \subset N$ compact such that $ \Gamma_N F_N =N$.
Then $F^o_N \not=\emptyset$. Define $\mathbb{F}_A=\{a(\mathbf{t})\mid \forall j\in \{1,\ldots ,r\}\, :\, |t_j|\le 1\}$
and
\[\mathbb{F}_Q:=F_A F_N B \subset Q\, . \]
Then $\mathbb{F}_{Q}$ is compact and $\mathbb{F}^o\not=\emptyset$. Furthermore
$\Gamma \mathbb{F}_Q=Q$ as $\Gamma_A\mathbb{F}_A=A$, $\Gamma_N\subset \{a\gamma a^{-1}\mid a\in
\mathbb{F}_A\, , \, \gamma\in \Gamma_N\}$, and $\Gamma_N\mathbb{F}_N=N$.

\begin{theorem}\label{thmain2} Assume that $Q=ANR$ is a closed subgroup
of $\GLR$ such that (1) -- (5) above holds.
Suppose that $\Gamma_N$ is a discrete subgroup of $N$ such that
$\Gamma_N\backslash N$ is compact.  Set $\Gamma=\Gamma_A
\Gamma_N$ and define
$\mathbb{F}_Q=\mathbb{F}_A\mathbb{F}_NB$ as above. Assume further that $\cO\subset \R^n$
is
an open  $Q$-orbit such that $\overline{\cO}\setminus \cO$ has
measure zero . Let $\omega\in \cO$ and assmume
that $Q^\omega=L \subset R$. Let $\mathbb{F}=\mathbb{F}_Q\cdot \omega$.
Then $(\Gamma ,\mathbb{F})$ is a frame generator.
\end{theorem}

\begin{proof} This follows from the discussion just before the statement of the Theorem.
\end{proof}

By condition (6) we see that the group $N$ as constructed in \cite{FO02} is
a simply connected
nilpotent Lie group. We recall here the most general statement about the existence
of co-compact subgroups of nilpotent Lie group, but first let us recall the
following definition. Let $\fg$ be a Lie algebra and let $\{X_1,\ldots ,X_r\}$ be
a basis for $\fg$. Then we can write
\[[X_i,X_j]=\sum_{k=1}^rc_{ijk}X_k\, .\]
The constants $c_{ijk}$ ($1\le i,j,k\le r$) are called the \textit{structure
constants} of $\fg$ relative to the basis $\{X_1,\ldots ,X_r\}$.

\begin{theorem}[Malcev, 1949]\label{Mal49} Suppose that $N$ is a
simply connected nilpotent Lie group with Lie algebra $\fn$.
Then $N$ contains a co-compact discrete subgroup if and only if
$\fn$ has a basis with rational structure constants.
\end{theorem}

The following is also well known and follows from Theorem \ref{Mal49}. See
also \cite{BH62}, p. 511, for proof.

\begin{theorem} Assume that $N\subset \GLR$ is unipotent (i.e., the Lie
algebra is nilpotent) and defined over $\mathbb{Q}$. Let
$N_{\mathbb{Z}}:=N\cap \GL (n,\mathbb{Z})$. Then
$N_{\mathbb{Z}}\backslash N$ is compact.
\end{theorem}

We notice also the following simple application of the ideas in the proof
of Theorem \ref{thmain1} and Theorem \ref{thmain2}

\begin{lemma} Let $H\subset \GLR$ be a closed subgroup. Assume that
$\cO$ is an open $H$-orbit under the twisted action such that
$H^\omega$ is compact for $\omega \in \cO$. Assume
that $\overline{\cO}\setminus \cO$ has measure zero and that
there exist a co-compact discrete subgroup $\Gamma\subset H$ such
that $\Gamma\cap H^\omega=\{e\}$. Then there exits a compact subset
$\mathbb{F}\subset \cO$ such that $(\Gamma, \mathbb{F})$ is
a frame generator.
\end{lemma}

\begin{proof} Let $F\subset H$ be a compact subset such that such that $\Gamma F=G$.

\end{proof}

\section{Action of reductive groups}

In this section we recall the construction from \cite{FO02} of the group $Q$.
This in particular gives us a simple criteria for the existence of
a co-compact discrete subgroup $\Gamma_N$ in $N$. We refer to
\cite{FO02} for proofs.

\begin{definition}
A closed subgroup $H\subset\mathrm{GL}(n,\mathbb{R})$ is called
\textbf{reductive} if there exist a $x\in\mathrm{GL}(n,\mathbb{R})$
such that $xHx^{-1}$ is invariant under transposition, $a\mapsto a^T$.
\end{definition}

We assume from now on that $H$ is reductive. For simplicity we can then
assume that $H^T=H$. In order to handle cases
where $\cO = H\cdot \omega $ do not necessarily have compact stabilizers, we will
assume that there exist an involution $\tau :H\rightarrow H$ such
that
\[
H_{o}^{\tau}\subset H^{\omega }=\left\{  h\in H\mid h^T ( \omega)=\omega \right\}  \subset
H^{\tau}
\]
where
\[
H^{\tau}=\left\{  h\in H\mid\tau(h)=h\right\}
\]
and the subscript ``$_{o}$'' indicates the connected component containing
the identity element. Notice that the involution $\tau$ may depend on
the open orbit. We can assume that $L$ is also
invariant under transposition. Then $\theta:h\mapsto(h^{T})^{-1}$ and
$\tau$ commute. Let $K=\mathrm{O} (n)\cap H=\left\{ k\in H\mid k^{T}
=k^{-1}\right\} $. Then $K$ is a maximal compact subgroup of $H$ and
$L\cap K$ is a maximal compact subgroup of $L$. Denote by $\mathfrak{h}$
the Lie algebra of $H$, i.e., $
\mathfrak{h}=\left\{  X\in M_{n}(\mathbb{R})\mid
\forall t\in\mathbb{R}:e^{tX}\in H \right\}$.
Then $\mathfrak{h}$ decomposes as
\begin{eqnarray*}
\mathfrak{h}  & =&\mathfrak{k}\oplus\mathfrak{s}\\
& =&\mathfrak{l}\oplus\mathfrak{q}\\
& =&\mathfrak{k}\cap\mathfrak{l}\oplus\mathfrak{k}\cap\mathfrak{q}\oplus\mathfrak{s}\cap
\mathfrak{l}\oplus\mathfrak{s}\cap\mathfrak{q}\,
\end{eqnarray*}
where
$$\mathfrak{s}=\left\{ X\in\mathfrak{h}\mid X^T=X\right\}$$
is the
subspace of symmetric matrices, and
$$\mathfrak{q}=\left\{ X\in\mathfrak{h}
\mid\tau(X)=-X\right\}\, . $$
Notice that
\[
[\mathfrak{l},\mathfrak{q}]\subset\mathfrak{q\quad}\text{and}\quad [
\mathfrak{k},\mathfrak{s}]\subset\mathfrak{s\,.}
\]
Even if we are not going to use it, we would
also like to recall the isomorphism $\fq\simeq T_{\omega}(\cO)$ given
in the following way. For $X\in \fq$ define a derivation $D_X$ by
\[D_X(f):=\frac{d\, }{dt}f(e^{tX}\cdot \omega)|_{t=0}\, .\]
Then $\fq \ni X\mapsto D_X\in T_{\omega}(\cO )$ is a linear isomorphism. As $L$ fixes $\omega$
it follows that $L$ acts on $T_{\omega}(\cO)$. Denote by $\ell_h$ the map $\ell_h(\eta )=
h\cdot \eta$. Then for $h\in L$ we have
\begin{eqnarray*}
(d\ell_h)_\omega (D_X)f&=&
\frac{d\, }{dt}f(he^{tX}\cdot \omega)|_{t=0}\\
&=&\frac{d\, }{dt}f(e^{t\mathrm{Ad}(h)(X)}\cdot \omega)|_{t=0}\\
&=&D_{\mathrm{Ad}(h)(X)}(f)\, .
\end{eqnarray*}
Hence, the action of $L$ corresponds to the natural action of $L$ on $\fq$. In particular
it follows then that the tangent bundle $T(\cO)$ can be described as
the vector bundle $T(\cO)=H\times_L\fq$.

Recall the linear maps $\mathrm{Ad}(a),\mathrm{ad}(X):\mathfrak{h}
\rightarrow\mathfrak{h}$, $a\in H$, $X\in\mathfrak{h}$, are given by
$\mathrm{ad}(X)Y=XY-YX=[X,Y]$ and $\mathrm{Ad}(a)Y=aYa^{-1}$. Let $\mathfrak{a}$
be a maximal commutative subspace of $\mathfrak{s}\cap\mathfrak{q}$; thus
$XY-YX=0$ for all $X,Y\in\mathfrak{a}$. Then the algebra $\mathrm{ad}
(\mathfrak{a})$ is also commutative. Let $(\cdot \mid \cdot)$ be the inner product
on $\mathfrak{h}$ given by $(X\mid Y): =\mathrm{Tr}(XY^T)$.
Then a simple calculation shows that
$ (\mathrm{ad}(X)Y\mid Z)
 =(Y\mid\mathrm{ad}(X^T)Z)$.
Thus $\mathrm{ad}(X)^T=\mathrm{ad}(X^T)$. In particular, if
$X\in\mathfrak{s}$ then $\mathrm{ad}(X)$ is symmetric. It follows that we
can diagonalize the action of $\mathfrak{a}$ on $\mathfrak{h}$. Specifically,
for $\alpha\in\fa^*$ set
\[
\mathfrak{h}^{\alpha}=\left\{  Y\in\mathfrak{h}\mid\forall X\in\mathfrak{a}\,:\mathrm{ad}
(X)Y=\alpha(X)Y\right\}  \,.
\]
Let $\Delta=\left\{\alpha\in\mathfrak{h}^*\mid\mathfrak{h}^{\alpha}
\not =\left\{0\right\}\right\}\setminus\left\{0\right\}$. Notice
that the set $\Delta$ is finite. Hence, there is a $X_r\in\mathfrak{a}$
such that $\alpha(X_r)\not =0$ for all $\alpha\in\Delta$. Let
$\Delta^{+}=\left\{\alpha\mid\alpha(X_r)>0\right\}$ and
Let
\[
\mathfrak{n}=\bigoplus_{\alpha\in\Delta^{+}}\mathfrak{h}^{\alpha}\, .
\]
Let $\mathfrak{m}_{1}=\left\{  X\in\mathfrak{h}\mid [\mathfrak{a},X]=\left\{
0\right\}  \right\}  $, and $\mathfrak{m}=\left\{  X\in\mathfrak{m}_{1}\mid\forall
Y\in\mathfrak{a}:(X\mid Y)=0\right\}  $. Then $\mathfrak{m}_{1}=\mathfrak{m}\oplus
\mathfrak{a}$. Furthermore $\mathfrak{m}$, $\mathfrak{n}$, and
$\mathfrak{p}=\mathfrak{m}\oplus\mathfrak{a}\oplus\mathfrak{n}$ are subalgebras of $\mathfrak{h}$
and
\begin{eqnarray*}
\mathfrak{h}  & =&\mathfrak{k}+\mathfrak{p}\\
&=& \mathfrak{l}+ \fp\, .
\end{eqnarray*}
Notice that this is not a direct sum in general because $\mathfrak{k}\cap\mathfrak{p}
=\mathfrak{k}\cap\mathfrak{m}$ and $\mathfrak{l}\cap\mathfrak{p} =\mathfrak{l}\cap\mathfrak
{m}$. Let $\mathfrak{m}_{2}$ be the algebra generated by
$\mathfrak{m}\cap\mathfrak{s}$, i.e., $\mathfrak{m}_{2}=[\mathfrak{m}\cap\mathfrak{s},
\mathfrak{m}\cap\mathfrak{s}]\oplus\mathfrak{m}\cap\mathfrak{s}$.
Then $\fm_2$ is an
ideal in $\mathfrak{m}$ and contained in $\mathfrak{m}
\cap\mathfrak{l}$.

Let
\[
\mathfrak{r}:=\left\{  X\in\mathfrak{m}\mid\forall Y\in\mathfrak{m}_{2}:(X,Y)=0\right\}
=\mathfrak{m}_{2}^{\perp}.
\]
Then $\mathfrak{r}$ is an ideal in $\mathfrak{m}$ and $\mathfrak{m}=\mathfrak{r}\oplus
\mathfrak{m}_{2}$. Let
$$N_{K}(\mathfrak{a})=\left\{  k\in K\mid\forall X\in
\mathfrak{a}:\mathrm{Ad}(a)X\in\mathfrak{a}\right\}  $$
and
$$N_{L\cap K}(\mathfrak{a}
)=N_{K}(\mathfrak{a})\cap L\, .$$
Finally let
$$M_{K}=Z_{K}(\mathfrak{a})=\left\{
k\in K\mid\forall X\in\mathfrak{a}:aXa^{-1}=X\right\}
$$
and
$$M_{L}=Z_{L\cap K}=L\cap M_K
(\mathfrak{a})\, .$$
Then
\begin{equation*}\label{Weyl}
W=N_{K}(\mathfrak{a})/M_{K}
\quad\mbox{and}\quad
W_{0}=N_{K\cap L}(\mathfrak{a})/M_{L}\subset W\, .
\end{equation*}
are finite groups
For $0\le j\le k=\# W/W_0$ choose $s_j\in N_K$ such that $s_0=e$ and
by obvious abuse of notation
\begin{equation*}\label{represent}
W=\dot\bigcup s_jW_{0}\quad\text{(disjoint union)\, .}
\end{equation*}

Let $P=\left\{ a\in H\mid\mathrm{Ad}(a)\mathfrak{p}=\mathfrak{p}
\right\} $, $A=\left\{ e^
{X}\mid X\in\mathfrak{a}\text{ }\right\} $ and $N=\left\{ e^{X}\mid
X\in\mathfrak{n}\right\} $. Then $P, A$, and $N$ are closed subgroups of $H$,
and $A,N\subset P$.
Let $M_{2}$ be the group generated by $\exp(\mathfrak{m}_{2})$, and $R_
{o}=\exp(\mathfrak{r})$. Then $F=\exp(i\mathfrak{a})\cap K\subset M_{K}$ is
finite and such that $R=FR_{o}$ is a group. Furthermore
\[
R\times M_{2}\times A\ni(r,m,a)\mapsto rma\in Z_{H}(A)
\]
is a diffeomorphism. Notice that by definition $F$ is central in
$Z_{H}(A)$ and $mFm^{-1}=F$ for all $m\in N_{K}(\mathfrak{a})$. Let
$M=RM_{2}$. Then $P=MAN$. Furthermore each element $p\in P$ has a unique
expression $p=man$ with $m\in M$, $a\in A$, $n\in N$. The final step
is now to define
\begin{equation}\label{defQ}
Q:=RAN=ANR\subset P\, .
\end{equation}
Then $Q$ is a closed subgroup of $H$ with Lie algebra $\mathfrak{q}=\mathfrak{r}
\oplus\mathfrak{a}\oplus\mathfrak{n}$. Notice that $\mathfrak{h}=\mathfrak{l}+
(\mathfrak{r}\oplus\mathfrak{a}\oplus\mathfrak{n})$ and
$\mathfrak{l}\cap(\mathfrak{r}\oplus\mathfrak{a}\oplus\mathfrak{n})
=\mathfrak{l}\cap\mathfrak{r}$.

By Theorem 4.3 in \cite{FO02} we have the following:

\begin{theorem} Choose $e=s_0,s_1,\ldots ,s_k\in W$ such that
$W$ is the disjoint union of the cosets $s_jW_0$. Then
\[\dot{\bigcup}_{j=0}^k Qs_jL\subset H\]
is open and dense. Furhermore there exist an analytic function
$\psi :H\to \mathbb{C}$ such
that
\[H\setminus \bigcup_{j=0}^k Qs_jL=\{h\in H\mid \psi (h)=0\}\, .\]
In particular $H\setminus \bigcup_{j=0}^k Qs_jL$
has measure zero.
\end{theorem}

\begin{theorem} Let the notation be as above. Let
$\cO\subset \R^n$ be an open orbit such
that $L=H^\omega$ is a symmetric subgroup. Let $\cO_j=Q\cdot (s_j\cdot \omega)$.
Then
\[\dot{\bigcup}_{j=0}^k\cO_j\subset\cO\]
is open and
\[\cO \setminus\dot{\bigcup}_{j=0}^k\cO_j\subset\cO\]
has measure zero.
Furthermore the stabilizer in $Q$ of $\tilde{\omega}$
in $\cO_j$ ($1\le j\le k$) is compact.
\end{theorem}

\begin{example}[$\GLR$ acting on symmetric matrices]\label{GLsym}
Let $V$ be the space $\mathrm{Sym}(n.\R)$ of symmetric $n\times n$ matrices.
Then $V\simeq \mathbb{R}^{n(n+1)/2}$
Under this identification the standard inner product on $\mathbb{R}^{n(n+1)/2}$ corresponds to
the inner product $(X,Y):=\mathrm{Tr}(XY)=\mathrm{Tr}(XY^T)$ on $V$.
Define an action of $H=\mathrm{GL}(n,\mathbb{R})$ on $V$ by
$$a(X)=(a^{-1})^TXa^{-1}\, .$$
Then
$$a\cdot X=a^\theta(X)=aXa^T\, .$$
Each symmetric matrix is up to conjugation determined by the signature and rank. The set $V_{\mathrm{reg}}=
\{X\in V\mid \det X\not=0\}$ is open and dense in $V$ and has measure zero. Furthermore each matrix in $V_{\mathrm{reg}}$ is
conjugate to one of the matrices $I(0)=I_{n}$, $I(n)=-I_{n}$ or
$$I(p)=\left(\begin{matrix}I_{n-p} & 0\cr 0 &I_{n-p}\end{matrix}\right)$$
where $1\le p\le n-1$. Denote the corresponding orbit by $\mathcal{O}_p$.
Notice that $X\mapsto -X$ defines a $\mathrm{GL}(n,\mathbb{R})$ isomorphism
$\mathcal{O}_p\simeq \mathcal{O}_{n-p}$. The group $\mathrm{O}(p,n-p)$ is by definition
given by
\[\mathrm{O}(p,n-p)=\{g\in \mathrm{GL}(n,\mathbb{R})\mid gI(p)g^T=I(p)\}\]
where we use the notation $\mathrm{O}(n)=\mathrm{O}(n,0)=\mathrm{O}(0,n)$.
Notice that $\mathrm{SO}(p,n-p)$ is compact if and only if $p=n$.
It follows that only the orbits $\mathcal{O}_0$ and $\mathcal{O}_n$ satisfy
the condition (W4). As before we let $\theta : \mathrm{GL}(n,\mathbb{R})\to \mathrm{GL}(n,\mathbb{R})$
be the involution $\theta (g)=(g^{-1})^T$. For  $j=0,\ldots ,n$ we define
$\tau_p : \mathrm{GL}(n,\mathbb{R})\to\mathrm{GL}(n,\mathbb{R})$
by
\begin{equation}\label{deftau}
\tau_p(g)=I(p)\theta (g)I(p)\, .
\end{equation}
Then
\[\tau_p(g)=g\iff gI(p)g^T=I(p)\]
and, hence,
\[H^{\tau_p}=\{g\in \mathrm{GL}(n,\mathbb{R})\mid \tau_p(g)=g\}=\mathrm{O}(p,n-p)=H^{I(p)}\]
and hence the orbit $\mathcal{O}_p\simeq \mathrm{GL}(n,\mathbb{R})/\mathrm{O}(p,n-p)$ is a symmetric space.

By abuse of notation we denote the derived involutions on $M(n,\R)$ by the same letters, i.e.,
\[\theta (X)=-X^T\, ,\qquad \tau_p(X)=-I(p)X^TI(p)\, .\]
Then $\theta (\exp X)=\exp (\theta(X))$ and $\tau_p(\exp X)=\exp (\tau_p(X))$.
Define as in the last section:
\begin{eqnarray*}
\fk &= &\{X\in M(n,\R)\mid \theta (X)=X\}=\fo (n)\, ,\\
\fs &= &\{X\in M(n,\R)\mid \theta (X)=-X\}=\Sy (n,\R)\, ,\\
\fh &= &\{X\in M(n,\R)\mid \tau_p(X)=X\}=\fo (p,n-p)\, ,\\
\fq&=&\{X\in M(n,\R)\mid \tau_p(X)=-X\}\, ,
\end{eqnarray*}
where we leave out the dependence of $\fh$ and $\fq$ on $p$ as that should be clear in each case.
Then
\begin{eqnarray*}
\glR&=&\fk\oplus \fs\\
&=&\fh\oplus \fq\\
&=&\fk\cap \fh\oplus \fk\cap \fq\oplus \fs\cap \fh\oplus \fs\cap \fq\, .
\end{eqnarray*}
In this case the abelian subalgebra $\fa$ is given by:
\[\fa=\{d(t_1,\ldots ,t_n)\mid t_1,\ldots ,t_n\in \mathbb{R}\}\]
where $d(d_1,\ldots ,t_n)=d(\mathbf{t})$ stands for the diagonal matrix with diagonal entries
$t_1,\ldots ,t_n$. In this case $\fa$ is maximal abelian in $\fs$ and in fact maximal
abelian in $\fh$. Hence, the centralizer of $\fa$ in $\fk$ and $\fl$ is trivial. Thus we
will have $P=Q$ in this case.

A simple calculation shows that
\[ [d(\mathbf{t}),E_{ij}]=(t_i-t_j)E_{ij}\, ,\]
i.e., the matrices $E_{ij}$ are the joint eigenvectors of $\{\ad (d(\mathbf{t}))
\mid d(\mathbf{t})\in \fa\}$, with eigenvalues $t_i-t_j$. Define $\alpha_{ij}:\fa\to \R$
by $\alpha_{ij}(d(\mathbf{t}))=t_i-t_j$. Then
$\Delta =\{\alpha_{ij}\mid 1\le i,j\le n,\, i\not= 0\}$. Let
$\Delta^+=\{\alpha_{ij}\mid 1\le i<j\le n\}$. Then $\fn$ is the
Lie algebra of upper triangular matrices with zero on the main
diagonal:
\[\fn=\bigoplus_{1\le i<j\le n}\R E_{ij}\, .\]
Furthermore $
A=\exp (\fa)=\{d(e^{t_1},\ldots ,e^{t_n})=\{d(\mathbf{t})\mid t_j>0\}$
and $N=\exp (\fn)$ is the group of upper triangular matrices with
one on the main diagonal. So in particular $AN$ is the group of upper triangular
matrices with positive diagonal elements. In this case
the group $M$ is given by $M=\{d(\mathbf{\epsilon})\mid \epsilon_j=\pm\}\simeq \{-1,1\}^n$
and $P=MAN$ is the group of regular upper triangular matrices.
The group $N_K(A)$ is the
finite group of elements with  only one coefficient non-zero in each column and
row, and that non-zero element is either $1$ or $-1$.
Let $\mathfrak{S}_n$ be the group of permutations of $\{1,\ldots ,n\}$. Then
$W$ acts on $\fa$ by
\[\sigma\cdot d(t_1,\ldots ,t_n)=d(t_{\sigma^{-1}(1)},\ldots ,t_{\sigma^{-1}(n)})\, .\]
Notice that this action can be realized as the conjugation by the orthogonal matrix
\[s_{ij}=(I_n-E_{ii}-E_{jj}) +E_{ij}-E_{ji}\in M^\prime \, .\]
It follows that $W\simeq \mathfrak{S}_n$.
The group $L\cap K=\mathrm{O}(p,n-p)\cap \mathrm{O}(n)$ is given by
\[L\cap K\simeq \mathrm{O}(p)\times \mathrm{O}(n-p)\]
where the isomorphism is given by
\[(a,b)\mapsto \left(\begin{matrix} a & 0\cr 0 & b\end{matrix}\right)\, . \]
Hence, $W_0\simeq \mathfrak{S}_p\times \mathfrak{S}_{n-p}$. In particular each
of the open $\GLR$-orbits is decomposed into $\frac{n!}{p!(n-p)!}$ $P$-orbits.
Finally we remark that in this example we can take
\[\Gamma_N=N\cap \GL (n,\mathbb{Z})\]
the group of upper triangular matrices with integer coefficients and
one on the diagonal.
\end{example}

\begin{remark}
Examples of pairs $(H,\mathbb{R}^n)$ such that $H$ is reductive and
has finitely many open oribts of full measure is given
by the \textit{pre-homogeneous vector spaces of parabolic type} (see
\cite{BR}). But there are simple examples where the stabilizer is not
symmetric. For that let  $H=\mathrm{SL}(2,\R)$. Then $H$ acts on
$\R^2$ in a natural way and the orbits are $\{0\}$ and
$\R^2\setminus \{0\}$. Hence there is only one open orbit, and
that orbit has full measure. The stabilizer of $e_1$
is the group
$$N=\left\{\left(\begin{matrix} 1 & x\cr 0 & 1\end{matrix}\right)\mid
x\in \R\right\}$$
which is not symmetric in $\mathrm{SL}(2,\R)$. Let
$$Q=\left\{\left(\begin{matrix} a & 0\cr y & 1/a\end{matrix}\right)\mid
a\not= 0\, y\in \R\right\}\, $$
Then $Q$ has three orbits $\{0\}$, $\{(0,y)^T\in \R^2\mid y\not=0\}$ and the open oribt
$\mathcal{O}=\{(x,y)^T\in \R^2\mid x\not= 0\}$. The stabilizer of $e_1\in \mathcal{O}$
is trivial and, hence, we can replace $\mathrm{SL}(2,\R)$ by $Q$ to construct
frames. Notice that $Q$ is isomorphic to the $(ax+b)$-group.
\end{remark}

\begin{remark} Assume that $(H,\mathbb{R}^n)$ is a pre-homogeneous vector space
of parabolic type. Then one can show that the same group $Q$ works
for all the open oribts \cite{BR}. Furthermore the group
$Q$ is admissible in the sense of  Laugesen,  Weaver,  Weiss, and  Wilson \cite{LWWW2002}
and contains an expansive matrix.
\end{remark}

\end{document}